\magnification 1200
\topskip=20pt
\vsize=19,5cm
\hsize=13cm

\def\Ref{{\rm {Ref}}} 

\def\card {{\rm {card}}}

 \def\Aut {{\rm {Aut}}}
 
 \def\S {{\Sigma}}

\def\Int {{\rm {Int}}}

\def\Z {\bf Z}

\def\mod {{\rm {mod}}}

 \def\dim  {{\rm {dim}}}

\def\mod {{\rm {mod}}\,}

\def\id {{\rm {id}}}

 \def\Ker {{\rm {Ker}}} \def\Im {{\rm {Im}}}

\def\Aut {{\rm {Aut}}}  
\def\H {{\rm {Homeo}}}  
\def\S {{\Sigma}}  
\def\W {{W}}  

\openup 1\jot

\headline{\ifnum\pageno=1 \hbox{}\else
\ifnum\pageno<10 \hbox to\hsize{\tenrm\hfil
\expno-0\folio \hfil}\else
\hbox to\hsize{\tenrm\hfil \expno-\folio \hfil}\fi\fi}
\nopagenumbers

\def\expno{878}
\noindent\hbox to\hsize{S\'eminaire BOURBAKI \hfil Juin 2000} \par
\noindent 52\`eme ann\'ee, 1999-00, n$^{\rm o}$ \expno \par
\vskip 1.9truecm

\centerline{\bf FAITHFUL LINEAR REPRESENTATIONS OF THE BRAID GROUPS} \par
\medskip
\centerline{by {\bf Vladimir TURAEV}}
\vskip 1.9truecm

The braid group on $n$ strings $B_n$  can be defined as the group
generated by $n-1$ generators
$\sigma_1,..., \sigma_{n-1}$ with defining relations
$$\sigma_i\sigma_j=\sigma_j\sigma_i \eqno (0.1)$$
if $\vert i-j\vert \geq 2$ and 
$$\sigma_i\sigma_{i+1} \sigma_i= \sigma_{i+1} \sigma_i\sigma_{i+1} \eqno (0.2)
$$
for $i=1,...,n-2$. This group  was introduced by  Emil Artin  in 1926. It
  has various    interpretations,   specifically,  as  the group of geometric braids in ${\bf R}^3$, as the mapping class group of
an
$n$-punctured   disc,  as the fundamental group of the configuration space of $n$ points on the plane,  etc. The algebraic properties of
$B_n$ have been studied by   many authors. To  mention  a few, 
 note  the solution of the conjugacy problem in $B_n$ given by  F. Garside, the papers of N. Ivanov 
and J. McCarthy who proved that the mapping class groups and in particular 
$B_n$ satisfy the
\lq\lq Tits alternative", and the  work of P. Dehornoy
establishing the existence of a left-invariant total order  on $B_n$ (see [Ga], [Iv2], [Ka]).

One of the most intriguing problems in the theory of braids is the
question of whether $B_n$ is 
linear, i.e., whether it admits a faithful representation into a group of matrices over a commutative ring.
This question  has its origins in a number of  interrelated facts and first of all in the discovery by  
W. Burau [Bu] of an $n$-dimensional linear representation of $
B_n$ which for a long time had been considered as a candidate for  a faithful representation.
  However, 
 as it was established by J.
Moody in 1991  this representation is not faithful for $n\geq 9$.  Later it was shown to be 
unfaithful for $n\geq 5$.  Thus, the question of the linearity of $B_n$
remained  open.

In 1999/2000 there appeared   a series of papers of D. Krammer and S. Bigelow who proved 
that   $B_n$ is linear for all $n$.  First there appeared a paper of   Krammer  [Kr1]  in which 
he constructed a homomorphism from $B_n$ to $GL( n(n-1)/2, R)$ where 
$R={{\bf Z}}[q^{\pm 1}, t^{\pm 1}]$ is the ring of Laurent polynomials on two variables.
He proved that this homomorphism  is injective for $n=4$ 
and conjectured that the same is true for all $n$. Soon after that,   
Bigelow   [Bi2]  gave a beautiful topological proof of this conjecture. Another proof based on  different 
ideas was   obtained by  Krammer [Kr2] independently.

The ring ${{\bf Z}}[q^{\pm 1}, t^{\pm 1}]$ can be embedded in the field of real numbers by assigning
to $q,t$ any algebraically independent non-zero real values. Therefore  $B_n$ embeds in 
$GL( n(n-1)/2, {{\bf R}})$. 
 As an   application, note that the linearity of $B_4$ implies the linearity of the group 
$\Aut (F_2)$, where $F_n$ is a free group of rank $n$ (see [DFG]). It is known that 
$\Aut (F_n)$
is not linear for $n\geq 3$, see [FP].

The representation of $B_n$ considered by  Krammer  and   Bigelow  is one of a   family of representations introduced
earlier by   R. Lawrence [La]. Her work was inspired by a study of the Jones   polynomial  of links and was concerned
with   representations of Hecke algebras arising from the actions of braids on homology of  configuration spaces.

The same representation of $B_n$ arises from a study of the so-called Birman-Murakami-Wenzl algebra $C_n$.
This algebra is a quotient of the group ring ${\bf C} [B_n]$ by certain relations inspired by the
theory of link polynomials. 
The irreducible finite dimensional representations of $C_n$ were described  by H. Wenzl in terms of Young diagrams. 
These representations yield  irreducible finite dimensional representations of $B_n$.
One of them was shown  by M. Zinno [Zi] and independently by V. Jones to be equivalent to the  
Krammer representation which is henceforth irreducible.

The aim of this paper is to present these results. In Sect.\  1 we consider the Burau representation and
explain why it is not faithful. In Sect.\   2 we   outline  Bigelow's approach following [Bi2].
In Sect.\   3 we  discuss the work of Krammer [Kr2]. Finally in Sect.\   4 we discuss the
Birman-Murakami-Wenzl algebras and the    work
of  Zinno. 

The author is indebted to S. Bigelow for useful comments on his work and to Ch. Kassel for careful
reading of a preliminary version of this paper.
 
\vskip 1.2truecm
\noindent {\bf 1. THE BURAU REPRESENTATION OF THE BRAID GROUP}
\vskip 0.5truecm
\par \noindent {\bf 1.1. Mapping class groups.}
It will be convenient for us   to view the braid group as the mapping class group 
of a punctured disc.     We  recall  here the definition  and a few
simple properties of  the mapping class groups.

Let $\S$ be a connected oriented surface.  By a {\it self-homeomorphism} of $\S$ we   mean an orientation preserving
homeomorphism $\S\to \S$ which fixes $\partial \S$ pointwise. Two such homeomorphisms are {\it  isotopic}  if they can be
included in a continuous one-parameter family of self-homeomorphisms of $\S$. The mapping class
group 
$\H(\S)$ of 
$\S$ is the group of isotopy classes of   self-homeomorphisms of $\S$ with the group operation   determined by
composition. 

 Each self-homeomorphism   of $\S$ induces an automorphism   of the abelian group
$H=H_1(\S; {{\bf Z}})$.  This is a    \lq\lq homological" representation $\H(\S)\to \Aut (H)$.
The action of homeomorphisms preserves the skew-symmetric bilinear  form $H
\times H \to
{{\bf Z}}$ determined by the
  algebraic intersection number. The value $[\alpha] \cdot [\beta] \in {{\bf Z}}$ of this form on  the homology classes $[\alpha], 
[\beta]
\in H$  represented by oriented
loops $\alpha,\beta$ on
$\S$  is computed as follows. 
Assume that $\alpha$ and $\beta$  lie in a generic position so that  they meet each
other  transversely in a finite set of points  
which are not self-crossings of $\alpha$ or $\beta$.
Then $$[\alpha] \cdot  [\beta]=\sum_{p\in \alpha\cap \beta} \varepsilon_p $$
where $\varepsilon_p=+1$ if the tangent vectors of $\alpha,\beta$ at $p$ form a positively oriented basis 
and
$\varepsilon_p=-1$ otherwise.

 An example of   a  self-homeomorphism of $\S$ is provided by the Dehn twist 
$\tau_\alpha$ about a simple closed curve
$\alpha\subset \S$.  It is defined as follows. Identify a regular neighborhood of $\alpha$ in $\S$ with the cylinder
$S^1\times [0,1]$ so that $ \alpha=S^1\times (1/2) $.  We   choose this identification so that the product of the
 counterclockwise orientation on 
$S^1=\{z\in {\bf C}\,\vert\, \vert z\vert =1\}$ and the right-handed orientation on $[0,1]$ corresponds to the given
orientation on $\S$.  The Dehn twist $\tau_\alpha: \S \to \S$ is the identity map outside  
$S^1\times [0,1]$ and sends any $(x,s) \in  S^1\times [0,1]$ to  $(e^{2\pi i s} x, s) \in  S^1\times [0,1]$.
To compute the action of $\tau_\alpha$ in homology, we orient $\alpha$  in an arbitrary way.  The effect of
$\tau_\alpha$ on an oriented   curve  transversal to $\alpha$ is   to insert   
$\alpha^{\pm 1}$  at each crossing  of $\alpha$ with this curve, where $\pm 1$ is the sign of
the crossing. Therefore   for any
$g\in H$,  
$$(\tau_\alpha)_* (g) =g+ ([\alpha] \cdot  g) [\alpha]   \eqno (1.1)$$
  Note that $\tau_\alpha$ and its action on $H$ do not depend on the choice of orientation on $\alpha$. 

A similar construction applies to   arcs in $\S$ whose endpoints are punctures. 
Assume that
  $\S$ is obtained from  another surface $\S'$ by puncturing, i.e., by removing a 
finite set of points lying 
in   $\Int (\S')$. These points will be called {\it punctures}. Let  $\alpha$ be an embedded arc
 in $\S'$
whose endpoints are  punctures $x_1,x_2$ and whose interior lies  in $\S$. 
One can define the Dehn \lq\lq half-twist"  $\tau_\alpha:\S\to \S$ which is the identity map
 outside a regular
neighborhood of $\alpha$ in $\S'$ and which exchanges $x_1$ and $x_2$.  This homeomorphism 
is obtained by  the isotopy of 
the identity map of $\S'$  rotating $\alpha$ about  its midpoint to the angle of 
$\pi$ 
in the direction provided by the orientation of  $\S$. Restricting the resulting 
 homeomorphism of $\S'\to \S'$
to $\S$ we obtain  $\tau_\alpha$.   To compute the action of $\tau_\alpha$ on   $ H=H_1(\S)$ we orient $\alpha$ from $x_1$
to $x_2$ and associate to    $\alpha$ a loop $\alpha'$ in $\S$ as follows.  Choose a point
$z\in
\alpha$ and  for   $i=1,2$  denote by
$\mu_i$   the loop in $\S$ beginning  at  $z$ and  moving
along 
$\alpha$ until coming very closely to $x_i$,  then encircling $x_i$ in the direction determined by
the orientation of $\S$ and finally  moving back   to $z$ along 
$\alpha$.    
Set  $\alpha'=\mu_1^{-1}
\mu_2$. The homotopy class  of  the loop $\alpha'$ on $\S$ does not
depend on the choice of
$z$. The effect of
$\tau_\alpha$ on an oriented   curve  transversal to $\alpha$ is  to insert   
$(\alpha')^{\pm 1}$  at each crossing   of $\alpha$ with this curve. Thus for any
$g\in H$, we have 
$$(\tau_\alpha)_* (g) =g+ ( [\alpha]\cdot g) [\alpha']  \eqno (1.2)
$$ where $[\alpha]\cdot g=-g\cdot [\alpha]\in {{\bf Z}}$ is the  algebraic intersection  number of 
$g$ with the 1-dimensional homology class $[\alpha] $ of $\S$ \lq \lq modulo infinity" represented by
$\alpha$.  

In general, the   action of   $\H(\S)$ on 
$H=H_1(\S )$  is  not
faithful.   We point out one   source  of non-faithfulness.
If $\alpha, \beta \subset \S$ are    simple closed curves with
$[\alpha] \cdot  [\beta]=0$
then formula  (1.1) implies that  $(\tau_\alpha)_*$ and $(\tau_\beta)_*$   commute
in $\Aut (H)$. The Dehn twists $\tau_\alpha, \tau_\beta$  themselves commute if and
only if   $\alpha$ is isotopic to a simple closed curve disjoint from
$\beta$, see  for instance [Iv1]. It is   easy to give examples of   simple closed curves 
$\alpha, \beta \subset \S$ which are not
disjoint up to isotopy but have zero algebraic intersection number.  
Then the commutator
$[\tau_\alpha,  \tau_\beta]$
lies in the kernel of the homological representation.    Using (1.2), one can  similarly derive elements of
the kernel    from  embedded arcs with endpoints in punctures.

\vskip 0.5truecm
\par \noindent {\bf 1.2. Braid groups.} Let $D=\{z\in {\bf C}\,\vert\, \vert z\vert \leq 1\}$ be the unit
disc with counterclockwise orientation. Fix  a set of $n \geq 1$ distinct punctures 
$X=\{x_1,...,x_n \} \subset  \Int (D)$. We shall assume that $x_1,...,x_n \in  (-1,+1)={\bf R}
\cap
\Int (D)$ and 
$x_1<x_2<...<x_n$.  Set $D_n=D\backslash X$. The group $\H (D_n)$  is denoted
$B_n$ and called the {\it $n$-th braid group}.  An element  of $  B_n$ is  an isotopy class
of  a homeomorphism $ D_n \to D_n$ which fixes $\partial D_n=\partial D=S^1$ pointwise.
Such a homeomorphism uniquely extends to a homeomorphism $D\to D$ permuting 
$x_1,...,x_n$. 
This defines a   group homomorphism from $B_n$
onto the symmetric group $S_n$.  We can equivalently define $B_n$ at the group of
isotopy classes of homeomorphisms $D\to D$ which fix $\partial D$ pointwise and preserve
$X$ as a set.

For $i=1,...,n-1$, the linear interval $[x_i,x_{i+1}] \subset (-1,+1) \subset  {\bf R} $ is   an embedded arc in
$D$ with endpoints in the punctures $x_i,x_{i+1}$. The corresponding 
Dehn  half-twist $D_n \to D_n$ is denoted by $\sigma_i$. 
It is a classical fact that  $B_n$ is generated by
$\sigma_1,..., \sigma_{n-1}$ with defining relations (0.1), (0.2). 
The image of $\sigma_i$ in $S_n$ is   the permutation $(i,i+1)$. 

Another  definition  of $B_n$ can be given in terms of    braids. A (geometric) {\it braid on $n$ strings} 
is an $n$-component  one-dimensional  manifold   $ E\subset D\times [0,1]$ such that 
 $E$ meets $D\times \{0,1\}$  orthogonally along the set $(X\times 0) \cup (X\times
1)$ and the projection on $[0,1]$
maps each component of $E$ homeomorphically onto $[0,1]$.  The  braids are
considered up to   isotopy in $D\times [0,1]$ constant on the endpoints. The group operation 
in the set of  braids is defined by glueing one braid on the top of the other one and
compressing the result into $D\times [0,1]$. 

The equivalence between these two 
definitions of $B_n$  is established as follows.  Any  
homeomorphism
$h:D\to D$ which fixes $\partial D$ pointwise is related to the identity map $\id_D$ 
by an isotopy $\{h_s:D\to D\}_{s\in [0,1]}$  such that $h_0=\id_D$ and $h_1=h$.
If $h(X)=X$ then the set $\cup_{s\in [0,1]} (h_s(X) \times s) \subset D\times [0,1]$
is a   braid. Its isotopy class depends only on the element of $B_n$
represented by $h$. This establishes an isomorphism between $B_n$ and the group of
  braids on $n$ strings.  The generator  $\sigma_i\in B_n$ corresponds to the $i$-th
\lq\lq elementary"  braid  represented by  a  plane diagram  consisting of $n$ linear intervals which are disjoint except at  one  
intersection point where  the $i$-th interval goes over the  $(i+1)$-th  interval.

\vskip 0.5truecm
\par \noindent {\bf 1.3. The Burau representation.}
Let $\Lambda$ denote the   ring  
${{\bf Z}}[t, t^{-1}]$. The Burau
representation  $B_n \to GL_n(\Lambda)$ sends the 
$i$-th   generator $\sigma_i\in B_n$ into the matrix
$$I_{i-1} \oplus  \pmatrix{  1-t &  t \cr
        1 &0 \cr}  \oplus I_{n-i-1}$$
where $I_k$ denotes the  identity $(k\times k)$-matrix and the non-trivial $(2\times 2)$-block       
appears in the $i$-th  and $(i+1)$-th  rows and columns (see [Bu]).   Substituting $t=1$, we obtain
the standard representation of the symmetric group $S_n$ by permutation matrices or equivalently the homological action of
$B_n$  on
$H_1(D_n)={\bf Z}^n$. The  Burau representation is  reducible: it splits as a direct sum of  an $(n-1)$-dimensional
representation and the  trivial one-dimensional representation.

The Burau representation is  known to be  faithful for $n\leq 3$, see [Bir].   J. Moody  [Mo] proved in 1991 that  it is not faithful 
for $n\geq 9$.   D. Long and M. Paton [LP]  extended Moody's argument  to $n\geq 6$.  Recently,
 S. Bigelow  [Bi1] 
proved that this representation is not faithful for $n=5$.  The case
$n=4$ remains open.

The geometric idea allowing to detect  non-trivial elements in the kernel of the Burau
representation is parallel to the one  at the end of  Sect.\   1.1. We first give 
a homological description  of the   Burau representation. Observe that $H_1(D\backslash  x_i ) ={\bf Z} $ is  
generated by the class of a small  loop encircling $x_i$ in the counterclockwise direction.
Each loop   in
$D\backslash
x_i$ represents $k $ times the generator where $k $ is the   {\it winding number} of  the loop 
around $x_i$. Consider the homomorphism $H_1(D_n) \to {{\bf Z}}$ sending
the homology class of a loop to  its {\it  total winding number} defined as  the sum of its winding numbers
around $x_1,...,x_n$. Let $\tilde D_n\to D_n$ be the  corresponding regular covering. The group of
covering transformations of $\tilde D_n$ is ${\bf Z}$ which we write as a multiplicative group 
with generator $t$. The  group $H_1( \tilde D_n)$  acquires thus the structure of a
$\Lambda$-module. It is easy to check that  this  is a free  
$\Lambda$-module of rank $n-1$. 

Fix a basepoint   $ d\in  \partial D$.  Any homeomorphism $h:D_n\to D_n$ representing an
element of $B_n$   lifts uniquely to a homeomorphism $\tilde h:\tilde D_n \to \tilde D_n$
which fixes the fiber over $d$ pointwise. This induces a $\Lambda$-linear automorphism 
$\tilde h_*$ of  $H_1( \tilde D_n)$. The map $h\mapsto \tilde h_*$  defines a  
representation   $B_n\to \Aut (H_1( \tilde D_n))$ equivalent to the $(n-1)$-dimensional Burau representation.

Now we   extend the algebraic intersection to arcs and  
refine
it so that it takes values in $\Lambda$.  Let $\alpha, \beta$ be two
embedded oriented arcs in $D_n$ with endpoints  in the punctures.  (We assume that  all four  
endpoints of $\alpha, \beta$ are distinct so that $n\geq 4$). Let
$\tilde
\alpha,
\tilde
\beta$ be   lifts of $\alpha, \beta$ to $\tilde D_n$, respectively. 
Set $$\langle \alpha, \beta \rangle= \sum_{k\in {\bf Z}} (t^k \tilde \alpha \cdot \tilde\beta) \,t^k \in \Lambda$$
where $t^k \tilde\alpha  \cdot \tilde \beta\in {\bf Z}$
is the   algebraic intersection number of the arcs $t^k \tilde \alpha, \tilde \beta$ in $\tilde D_n$. 
This finite sum  is only defined up to multiplication by  a power of $t$ depending on
the choice of $\tilde \alpha, \tilde \beta$. This will not be  important for us since
we  are only interested in whether or not  $\langle \alpha, \beta \rangle= 0$.
To compute 
$\langle \alpha, \beta \rangle$ explicitly one deforms $\alpha$ in general position with respect to $\beta$. Then 
$\langle
\alpha, \beta \rangle=\sum_{p\in
\alpha\cap
\beta}  \varepsilon_p \,t^{k_p}$   where $\varepsilon_p=\pm  $ is the  intersection sign  at
$p$ and $k_p\in {\bf Z}$. The exponents $\{k_p\}_{p }$
are determined by the following condition: if $p,q\in \alpha\cap \beta$,  
then $k_p-k_q$ is the total winding number of the loop going from $p$ to $q$ along
$\alpha$ and then from $q$  to $p$ along $\beta$.  Note that
$$\langle  \beta , \alpha  \rangle=
\sum_{k\in {\bf Z}} (t^k \tilde \beta \cdot \tilde\alpha) \,t^k
=\sum_{k\in {\bf Z}} (  \tilde \beta \cdot t^{-k} \tilde\alpha) \,t^k
=-\sum_{k\in {\bf Z}} (t^{-k} \tilde \alpha \cdot \tilde\beta) \,t^k= -\overline {\langle \alpha, \beta \rangle}$$
where the overline  denotes the   involution in $\Lambda$ sending any $t^k$ to $t^{-k}$. 
Hence $\langle \alpha, \beta \rangle=0$ if and only if $\langle  \beta , \alpha \rangle=0$.

We claim that if $\langle \alpha, \beta \rangle=0$, then the automorphisms
$(\tilde \tau_\alpha)_*, (\tilde \tau_\beta)_*$ of  $H_1(\tilde D_n)$ commute. 
Observe that the   loops $\alpha', \beta'$ on $D_n$ associated to $\alpha, \beta$ as in Sect.\   1.1 have zero total winding
numbers and therefore lift to  certain loops
$\tilde
\alpha',
\tilde \beta'$ in  $\tilde D_n$.  The effect of $\tilde \tau_\alpha$ on any  oriented loop
$\gamma$ in $\tilde D_n$ is to insert a lift of $(\alpha')^{\pm 1}$ at each crossing of $\gamma$ with
the preimage of $\alpha$ in $\tilde D_n$.   Thus
$$(\tilde \tau_\alpha)_* ([\gamma]) = [\gamma] + \lambda_{\gamma} \, [\tilde \alpha']  $$ for a certain Laurent polynomial  $\lambda_{\gamma} \in
\Lambda$.  The coefficients of  $\lambda_{\gamma}$ are the algebraic intersection numbers of $\gamma$
with   lifts of $\alpha$ to $\tilde D_n$. By
$\langle \alpha, \beta \rangle=0$, any lift of $\alpha$ has algebraic intersection number zero with any
lift of $\beta$ and hence with any lift of $\beta'$.  Therefore, 
$\lambda_{\tilde \beta'}=0$ and $(\tilde \tau_\alpha)_* ([\tilde \beta'])=[\tilde \beta']$. Similarly, 
$(\tilde \tau_\beta)_* ([\gamma]) = [\gamma] + \mu_{\gamma}  \, [\tilde \beta']  $  with  $\mu_{\gamma}  \in \Lambda$ and 
$(\tilde \tau_\beta)_* ([\tilde \alpha'])=[\tilde \alpha']$.
We conclude that
$$ (\tilde \tau_\alpha \tilde \tau_\beta)_* ([\gamma])
=[\gamma] + \lambda_{\gamma} \, [\tilde \alpha']   + \mu_{\gamma} \, [\tilde \beta']
= (\tilde \tau_\beta \tilde \tau_\alpha)_* ([\gamma]).$$

To show   that the  Burau representation is not faithful it remains to
  provide an example of    oriented  embedded 
  arcs $\alpha, \beta$ in $D_n$ with endpoints  in distinct punctures  such that $\langle \alpha, \beta
\rangle=0$ and $\tau_\alpha   \tau_\beta \neq  \tau_\beta
 \tau_\alpha$ in $B_n$.  
 For   $n\geq 6$, the simplest known  example  (see  [Bi1]) is provided by the pair  
$\alpha=\varphi_1 ([x_3, x_4] ),\, \beta=\varphi_2 ([x_3, x_4] ) $ where
$$\varphi_1=\sigma_1^2  \sigma_{2}^{-1} \sigma_{5}^{-2} \sigma_4 , \,\,\,\,\,\, \varphi_2=   \sigma_{1}^{-1} 
\sigma_2
 \sigma_5 \sigma_{4}^{-1}.$$  To compute $\langle \alpha, \beta \rangle $ one can draw 
$\alpha, \beta$ and use  the recipe above. To prove that the braids $\tau_\alpha = \varphi_1 \sigma_3
\varphi_1^{-1}$ and 
$ \tau_\beta=  \varphi_2 \sigma_3 \varphi_2^{-1}$ do not commute,  one can use 
the solution of the word problem in $B_n$ or 
  the methods of the Thurston theory of surfaces (cf. Sect.\  2.2). Thus the commutator
$[\tau_\alpha,  \tau_\beta]$
lies in the kernel of the Burau representation.  
   This commutator can be represented by a word of length 44 in the  
generators   $\sigma_1,..., \sigma_{5}$.

Similar ideas apply in the case $n=5$, although one  has to extend them to arcs  relating the punctures
to    the base point  
$d\in \partial D$.   The  shortest known word in the generators $\sigma_1,..., \sigma_{4}$ representing an element of
the kernel has length 120.

\vskip 1.2truecm
\noindent {\bf 2. THE WORK OF  BIGELOW}
\vskip 0.5truecm
\par \noindent {\bf 2.1. A  representation of $B_n$.} We use the notation $D, D_n, X=\{x_1,...,x_n\}$  introduced in
Sect.\   1. Let   $C$ be the space of all    unordered  pairs of distinct points in $D_n$.
This space is obtained from    $(D_n\times D_n)\backslash \hbox{diagonal}$ by the 
identification $\{x,y\}= \{y,x\}$ for  any  distinct $x, y \in D_n$. It is clear that $C$ is a  connected
non-compact 4-manifold  with boundary. 
It has  a  natural orientation induced by the counterclockwise  orientation of
$D_n$.  Set $d=-i\in \partial D$ where  $i=\sqrt{-1} $ and  $d'= -i\, e^{  { {\varepsilon \pi i} \over {2}}}\in \partial D$
  with small
positive $\varepsilon $. We take
$c_0=\{d,d'\}$ as the basepoint for $C$.

A  closed curve $\alpha:[0,1] \to C$   can be written in the form
$\alpha(s)= \{\alpha_1(s),\alpha_2(s)\}$ where $s\in [0,1]$ and $\alpha_1,\alpha_2$ are arcs in $D_n$ such that
$\{\alpha_1(0),\alpha_2(0)\}=\{\alpha_1(1),\alpha_2(1)\}$.
The   arcs  $\alpha_1,\alpha_2 $ are either  both loops or can be composed with each other. They form thus a closed
oriented one-manifold   mapped to 
$D_n$.  Let $a(\alpha)\in {\bf Z}$ be the total winding number of this one-manifold around  the punctures  
$\{x_1,...,x_n\}$.  Composing the map $s\mapsto  (\alpha_1(s)-\alpha_2(s))/ \vert \alpha_1(s)-\alpha_2(s)\vert: [0,1] \to
S^1$ with the projection $S^1\to {\bf R} {\bf P}^1 $ we obtain a loop in $RP^1$.
The corresponding element of $H_1 (RP^1)={\bf Z}$ is denoted by $b(\alpha)$.
The formula  $\alpha \mapsto q^{a(\alpha)} t^{b(\alpha)}$ defines 
a homomorphism, $\phi$, from $H_1(C)$ to the  (multiplicatively written) free abelian group with basis
$q,t$.   Let  $R={{\bf Z}}[q^{\pm 1}, t^{\pm 1}]$ be the group ring
of this group.

Let
$ 
\tilde C
\to C$ be a  regular covering corresponding to the kernel of   $\phi$.  The generators 
$q,t$ act on $\tilde C$  as commuting covering transformations. The homology group $H_2(\tilde C)=H_2(\tilde C;{\bf Z})$
becomes thus an
$R$-module. 

Any self-homeomorphism $h$ of $ D_n$  induces by $h(\{x,y\})=  \{h(x), h(y)\}$ a  homeomorphism   $C\to C$
also denoted    $h$. It is easy to check that $h(c_0)=c_0$ and the action   of $h$ on $H_1(C)$ commutes with $\phi$. 
Therefore this  homeomorphism  $C\to C$ lifts uniquely to a map $\tilde h: \tilde C \to \tilde C$ which fixes 
the fiber over $c_0$ pointwise and commutes with the covering transformations.  Consider the representation  
$B_n\to
\Aut  (H_2(\tilde C))$ sending the isotopy class of $h$ to the $R$-linear automorphism $\tilde h_*$ of $H_2(\tilde C)$.

\vskip 0.5truecm
\par \noindent {\bf 2.2. Theorem.} (S. Bigelow   [Bi2]) -- {\it  The representation $B_n\to \Aut (H_2(\tilde
C))$ is faithful for all
$n\geq 1$.}

\vskip 0.5truecm

We outline below the main ideas of Bigelow's proof. The proof uses almost no information about the structure
of the $R$-module $H_2(\tilde C)$. The only thing  needed is the absence of $R$-torsion or more precisely the fact that
multiplication by a  non-zero polynomial  of type $q^at^b-1$ has zero kernel in $H_2(\tilde C)$. In fact, $H_2(\tilde C)$ is a free
$R$-module of rank
$n(n-1)/2$, as it was essentially shown in  [La].

We shall use one well-known fact concerning isotopies of   arcs on surfaces. 
 Let $N, T$ be    embedded arcs in $D_n$ with distinct endpoints  lying either  in the punctures or on
$\partial D_n$.  Assume that the interiors of $N,T$ do not meet $\partial D_n$,     and that $N$ intersects
$T$  transversely  (in a finite number of points).  A
{\it  bigon} for the pair
$(N,T)$ is  an embedded disc in $\Int (D_n)$ whose boundary is formed by  one subarc of $N$ and one subarc of $T$
and whose interior is disjoint from $N$ and $T$. It is clear that in the presence of a bigon there is an isotopy of $T$ 
constant on the  endpoints and decreasing
$\#(N\cap T)$ by two. Thurston's theory of surfaces implies the converse: if there is an isotopy of $T$ (rel endpoints)
decreasing 
$\#(N\cap T)$ then the pair $(N,T)$ has at least one bigon, cf.  [FLP, Prop. 3.10].

\vskip 0.5truecm
\par \noindent {\bf 2.3. Noodles and forks.}  We need the following notation. For arcs  $\alpha, \beta:[0,1] \to D_n$
such that $\alpha(s)\neq \beta(s)$ for all $s\in [0,1]$,   we denote by  $\{\alpha,\beta\}$   the arc in $C$ given by
$\{\alpha,\beta\} (s)= \{\alpha(s), \beta(s)\}$. We fix once and forever a point  $\tilde c_0 \in \tilde  C$ lying over $c_0\in C$.

A {\it noodle}  in $D_n$  is
an embedded arc
$N\subset D_n$ with endpoints
$d$ and $d'$. For a noodle $N$,  the set   $\S_N=\{\{x,y\}\in C \,\vert \, x, y\in N, x\neq y\}$  
is  a surface in $C$ containing $c_0$.  
This surface is homeomorphic to a triangle with one edge removed.  We   orient  $N$ from $d$
to
$d'$ and orient 
$
\S_N$ as follows: at a point
$\{x,y\}=\{y,x\}\in \S_N$ such that $x$ is closer to $d$ along $N$ than $y$, the orientation of $ \S_N$ is the product  of the 
orientations of
$N$ at $x$ and  $y$ in this order. Let  $\tilde \S_N$  be the lift of $\S_N$ to $\tilde C$ containing $\tilde c_0$. The orientation  
of $\S_N$ lifts to $
\tilde
\S_N$ in the obvious way. Clearly,  $\tilde \S_N$ is a proper surface in $\tilde C$ in the sense that
$\tilde \S_N\cap \partial \tilde C=\partial \tilde \S_N$.

 A {\it
fork} in $D_n$ is an embedded tree
$F\subset D$ formed by three edges and  four vertices
$d,x_i,x_j, z$
such that $F\cap \partial D=d , F\cap X= \{x_i,x_j\}$ and  $z$ is a common vertex of all 3 edges. The edge, $H$, relating 
$d $ to $z$ is called the {\it handle} of $F$. The union, $T$, of the other two edges is an embedded arc with  endpoints
$ \{x_i,x_j\}$. This arc is called the {\it tines} of $F$. Note that in a small neighborhood of $z$, the handle $H$ lies on one side
of
$T$ which distinguishes  a  side of $T$. We    orient $T$ so that its distinguished side lies on its right. The handle $H$
also has a distinguished side    determined by
$d'$.
 Pushing  slightly  the graph
$F=T\cup H$ to the  distinguished side (fixing the vertices $  x_i,x_j $ and pushing $d$  to $d'$)   
 we obtain a \lq\lq parallel copy" $F'=T'\cup
H'$. The graph
$F'$ is a fork with handle
$H'$, tines
$T'$, and vertices $d',x_i,x_j, z'$ where $z'=T'\cap H'$  
lies on the distinguished side  of both
$T$ and
$H$. We   can assume that $F'$ meets $F$ only in common vertices $ \{x_i,x_j\}=T\cap T'$  and
in one point lying on
$H\cap T'$ close  to  $z,z'$. 
 
For a fork  $F$,  the set   $\S_F=\{\{y,y'\} \in C \,\vert \, y \in T \backslash \{x_i,x_j\},   y'  \in T' \backslash \{x_i,x_j\} \}$  
is  a surface in $C$ homeomorphic to $(0,1)^2$.  Let $\alpha_0 $  be an arc  from $d$ to $z$ along $H$
and let $\alpha'_0 $  be an arc  from $d'$ to $z'$ along $H'$. Consider  the arc
$\{\alpha_0, \alpha'_0\}$ in $ C$ and denote by  $ \tilde  \alpha$ its  lift to $\tilde C$ which starts in $\tilde c_0$. 
Let  $\tilde \S_F$ be the lift of $\S_F$ to $\tilde C$ which contains the lift $\tilde  \alpha (1)$ of the point $\{z,z'\}\in \S_F$.
 The
surfaces
$
\S_F$ and
$\tilde
\S_F$  have a  natural orientation determined by the orientation of
$T$ and  the induced  orientation of $T'$. 

We shall use the surfaces $\tilde \S_N, \tilde \S_F$ to  establish a duality between   noodles and forks. More precisely, for
any    noodle
$N$ and any fork $F$  we define an element  $\langle N,F \rangle $ of $R$ as follows.
   By applying a preliminary isotopy we can assume that $N$ intersects $T$ transversely in
$m\geq 0$ points
$z_1,...,z_m$ (the numeration is arbitrary; the intersection of $N$ with $H$ may be not transversal).  We  choose the parallel
fork
$F'=T'\cup H'$  as above   so that $T'$ meets $N$ transversely in $m$ points $z'_1,...,z'_m$ where each pair $z_i, z'_i$ is
joined by a short arc  in $N$ which  lies in the narrow strip bounded by $T\cup T'$ and meets no other $z_j, z'_j$. Then the
surfaces 
$\S_F$ and $ \S_N$ intersect transversely in $m^2$ points  $  \{z_i, z'_j\}$ where $ i,j=1,...,m$.
Therefore for any $a,b\in {\bf Z}$,
the image $q^a t^b \tilde \S_N$ of $  \tilde \S_N$ under the covering transformation $q^a t^b  $ meets
$\tilde \S_F$ transversely. Consider the algebraic intersection number $q^a t^b \tilde \S_N \cdot \tilde \S_F \in {\bf Z}$  
 and set 
$$\langle N,F \rangle = \sum_{a,b\in {\bf Z}} (q^a t^b \tilde \S_N \cdot \tilde \S_F)\, q^a t^b . $$
The sum on the right-hand side is finite (it has $\leq m^2$  terms) and  thus defines  an element of $R$. 

\vskip 0.5truecm
\par \noindent {\bf 2.4. Lemma.} -- {\it  $\langle N,F \rangle $ is invariant under isotopies of $N$ and $F$  in $D_n$ 
constant on   the endpoints.}

\vskip 0.5truecm
\par \noindent {Proof.} -- We first compute  $\langle N,F \rangle $ explicitly. Let
$N\cap T= \{z_1,...,z_m\}$ and $  N\cap T'= \{z'_1,...,z'_m\}$   as above. For every pair $i,j\in \{1,...,m\} $, 
there exist unique integers $a_{i,j}, b_{i,j} \in {\bf Z}$ such that 
$ q^{a_{i,j}} t^{b_{i,j}} \tilde \S_N $ intersects  $\tilde \S_F$ at a point   lying over 
$\{z_i,z'_j\}\in C$. Let $\varepsilon_{i,j}=\pm 1$ be the sign of that intersection. Then
$$\langle N,F \rangle=\sum_{i=1}^m \sum_{j=1}^m \, \varepsilon_{i,j} \,q^{a_{i,j}} t^{b_{i,j}}.\eqno  (2.1)$$
 
The numbers $a_{i,j}, b_{i,j}$ can be computed  as follows.
Let $\alpha_0  $  be an arc  from $d$ to $z$ along $H$
and let $\alpha'_0 $  be an arc  from $d'$ to $z'$ along $H'$.
Let $\beta_i $  be an arc  from $z$ to $z_i$ along $T$
and let $\beta'_j $  be an arc  from $z'$ to $z'_j$ along $T'$. Finally, let $ \gamma_{i,j}$ and $ \gamma'_{i,j}$ be disjoint
arcs in
$N$ connecting the points $z_i, z'_j$ to the endpoints of $N$. Note that 
$\delta_{i,j}=\{\alpha_0, \alpha'_0\} \{\beta_i, \beta'_j\} \{\gamma_{i,j}, \gamma'_{i,j}\}$ is a loop in $C$. Then 
$$q^{a_{i,j}} t^{b_{i,j}}=\phi (\delta_{i,j}).\eqno  (2.2)$$
Indeed, we can lift   $\delta_{i,j}$ to a path $\tilde  \alpha  \tilde \beta  \tilde \gamma $
 in $\tilde C$ beginning at  $\tilde c_0$ 
where $\tilde  \alpha , \tilde \beta , \tilde \gamma $ are lifts of $\{\alpha_0, \alpha'_0\}, \{\beta_i, \beta'_j\} ,\{\gamma_{i,j},
\gamma'_{i,j}\}$, respectively. By definition of $\tilde \S_F$, the point $\tilde  \alpha(1)=\tilde \beta (0)$ lies on $\tilde \S_F$. 
Hence the lift  $\tilde  \beta $ of  $ \{\beta_i, \beta'_j\}$ lies  on $\tilde \S_F$. The path $\tilde  \alpha \tilde
\beta \tilde \gamma $ ends at  
$ \phi (\delta_{i,j}) (\tilde c_0)=\tilde \gamma (1)$. Hence the 
lift  $\tilde  \gamma$ of 
$\{\gamma_{i,j}, \gamma'_{i,j}\}$  lies on   $ \phi (\delta_{i,j}) \tilde \S_N$. Therefore the point 
$\tilde  \beta (1)=\tilde  \gamma (0) $  lying over $\{z_i,z'_j\}$   belongs to 
  both  $\tilde \S_F$ and   $ \phi (\delta_{i,j}) \tilde \S_N$. This yields  (2.2).

 Note that the  residue  $b_{i,j} \,(\mod 2)$ depends on whether the two
points of $D_n$ forming a point in $C$ switch places moving along the loop $\delta_{i,j}$. This is determined by which of
$z_i, z'_j$ lies closer to
$d$ along
$N$.

To compute $\varepsilon_{i,j}$  we  observe that 
$\varepsilon_{i,j}$ is determined by the signs of the intersections of $N$ and $T$ at $z_i, z_j
$ and by which of $z_i,z'_j$ lies closer to $d$ on $N$. The sign of the intersection  of $N$ and $T$ at $z_i$
is  $+$ if $N$ crosses $T$    from  the left to the right and is $-$ otherwise.   By our choice of
orientations on $N$  and $T$, this  sign  is $+$ if $z_i$ lies closer to
$d$ along
$N$ than $z'_i$ and is $-$ otherwise. Hence, this sign is $(-1)^{b_{i,i}}$.Therefore 
$\varepsilon_{i,j}$ is determined by  $b_{i,i} +b_{j,j} +b_{i,j} \,(\mod 2)$.  A  precise   
computation shows that
$$\varepsilon_{i,j}= -(-1)^{b_{i,i}  + b_{j,j}+b_{i,j}}. \eqno (2.3)$$

Now we can prove the lemma. It suffices to  fix $F$ and to prove that 
$\langle N,F \rangle $ is invariant under isotopies of $N$. A  generic  isotopy of
$N$ in $D_n$ can be split into a finite sequence of local moves of   two kinds:   (i)  isotopies keeping $N$
transversal to
$T$, (ii) a move pushing a  small subarc of $N$ across a subarc of $T\backslash z$. It is clear from the discussion above that
the move (i)    does not change $\langle N,F \rangle $. The move (ii)   adds  two  new intersection points $z_{m+1}, z_{m+2}$ to
the set 
$N\cap T= \{z_1,...,z_m\}$.
It follows  from definitions  and the discussion above that for any $i=1,...,m+2$,
$$a_{i,m+1}=a_{i,m+2}, \,\,\,b_{i,m+1}=b_{i,m+2},  \,\,\, \varepsilon_{i,m+1}=-\varepsilon_{i,m+2}.$$ Hence 
for any $i=1,...,m+2$, the terms 
$\varepsilon_{i,m+1} q^{a_{i,m+1}} t^{b_{i,m+1}}$ and $\varepsilon_{i,m+2} q^{a_{i,m+2}} t^{b_{i,m+2}}$ cancel each other.
Similarly, for any $j=1,...,m$, the terms 
$\varepsilon_{i,j} q^{a_{i,j}} t^{b_{i,j}}$ with $i=m+1, m+2$  cancel each other.
Therefore   $\langle N,F \rangle $ is the same before and after the move.

\vskip 0.5truecm
\par \noindent {\bf 2.5. Lemma.} -- {\it   $\langle N,F \rangle =0$ if and only if   there is an isotopy  $\{T(s)\}_{s\in [0,1]}$ of
the tines $T=T(0)$ of $F$  in
$D_n$ (rel endpoints) such that   $T(1)$ is disjoint from $N$.}

\vskip 0.5truecm
\par \noindent {Proof.} -- Any isotopy $\{T(s)\}_{s\in [0,1]}$ of
  $T=T(0)$ extends to an ambient isotopy of $D_n$ constant on $\partial D_n$ and  therefore extends to an isotopy 
$\{F(s)\}_{s\in [0,1]}$ of
the fork $F=F(0)$. If $T(1)$ is disjoint from $N$ then by Lemma  2.4, 
$\langle N,F \rangle=\langle N,F(1)  \rangle=0$.

The hard part of the lemma is   the opposite implication. By applying a preliminary
isotopy, we  can assume that $T$ intersects $N$ transversely at a {\it  minimal} number of points $z_1,...,z_m$ with $m\geq 0$.
We assume that $m\geq 1$ and show that
$\langle N,F
\rangle \neq 0$. To this end we use the lexicographic ordering on monomials $q^a t^b$. Namely we write
$q^a t^b \geq q^{a'} t^{b'}$ with $a,b,a',b'\in {\bf Z}$ if either $a>a'$ or $a=a'$ and $b\geq b'$.  We say that the ordered pair
$(i,j)$ with
$i,j\in \{1,...,m\}$ is  {\it maximal} if $q^{a_{i,j}} t^{b_{i,j}}  \geq q^{a_{k,l}} t^{b_{k,l}}
$ for any
$k,l \in \{1,...,m\}$. We claim that

\vskip 0.5truecm

$(\ast)$  {\it if the pair $(i,j)$ is maximal, then   $b_{i,i}=b_{j,j}=b_{i,j}$.}

\vskip 0.5truecm

This  claim and  (2.3) imply that all  entries of the maximal monomial, say $q^a t^b$,  in (2.1) occur with
the same sign $-(-1)^b$.  Hence $\langle N,F
\rangle \neq 0$.

To prove $(\ast)$ we first  compute $a_{i,j}$ for any $i,j$ (not necessarily maximal). Let $\xi_i$ be the loop obtained  by
moving from $d$ to $z_i$ along $F$ then back to $d$ along $N$. Let $a_i$ be the total winding number of $\xi_i$ around all
$n$ punctures.
Let $\xi$ be the loop obtained by moving from $d$ to $d'$ along $N$, and  then  moving  clockwise along
$\partial D$ back to $d$. Let $a$ be the total winding number of $\xi$. We claim that
$$a_{i,j}=a_i+a_j+a.\eqno  (2.4)$$ 
Indeed, if $b_{i,j}$ is even then the  paths $\alpha_0 \beta_i \gamma_{i,j}$ and
$\alpha'_0\beta'_j\gamma'_{i,j}$  (in the notation of Lemma 2.4) are loops and 
$a_{i,j}$ is the sum of their total winding numbers. These loops are homotopic  in $D_n$ to
$\xi_i$ and $\xi_j \xi$, respectively. This implies (2.4). If $b_{i,j}$ is  odd then 
$a_{i,j}$ is the   total winding number of the loop 
$\alpha_0 \beta_i \gamma_{i,j} \alpha'_0\beta'_j\gamma'_{i,j}$.  This loop is homotopic in $D_n$ to $\xi_i  \xi  \xi_j$ which  
implies (2.4) in this case.

Suppose now that the pair $(i,j)$ is maximal. Then $a_{i,j}$ is maximal among all the integers $a_{k,l}$.
By (2.4),  it follows that $a_i =a_j$ is maximal among all the integers $a_k$. (Although we shall not
need it, observe that then
$a_{i,i}=a_{j,j}=a_{i,j}$).

We now show that $b_{i,i}=b_{i,j}$. By the  maximality of   $(i,j)$, we have that 
$b_{i,i}\leq b_{i,j}$.  Suppose, seeking a contradiction, that $b_{i,i}<b_{i,j}$. Let $\alpha$ be an embedded arc from $z'_i$ to
$z'_j$ along $T'$. Let $\beta$ be an embedded arc from $z'_j$ to
$z'_i$ along $N$. 

If $\beta$ does not pass through the point $z_i$, then we denote by $w$   the winding number of the loop
$\alpha \beta$ around
$z_i$. Observe that $b_{i,j}-b_{i,i}=2w$. To see this, consider   the loop 
$
\delta_{i,j}=\{\alpha_0, \alpha'_0\} \{\beta_i, \beta'_j\} \{\gamma_{i,j}, \gamma'_{i,j}\}$
  appearing in (2.2).  
Clearly,  $\beta'_j \sim \beta'_i \alpha$, where $\sim$ denotes   homotopy of paths in $D_n \backslash
z_i$ constant on the endpoints. The assumption that
$\beta$ does not pass through
$z_i$ implies that  $\gamma_{i,j}=\gamma_{i,i}$ and  $\gamma'_{i,j}\sim \beta  \gamma'_{i,i}$. 
Then $$\delta_{i,j}=\{\alpha_0, \alpha'_0\} \{\beta_i, \beta'_j\} \{\gamma_{i,j}, \gamma'_{i,j}\}\,\sim\,
\{\alpha_0, \alpha'_0\} \{\beta_i,
\beta'_i\} \{z_i, \alpha \beta \}
\{\gamma_{i,i},
\gamma'_{i,i}\}$$ where $z_i$ denotes the constant path in $z_i$. This implies that $b_{i,j}-b_{i,i}=2w$.

If   $\beta$   passes through   $z_i$, we  first modify $\beta$ in a small neighborhood of $z_i$ so that $z_i$ lies to its left.
Let
$w$ be the winding number of the loop 
$\alpha
\beta$ around
$z_i$. A  little more difficult but similar  argument shows that  $b_{i,j}-b_{i,i}=2w-1$. In either case
$w>0$. 

Let $D_0=D\backslash z_i$  and $p:\hat  D_0\to D_0$  be the universal (infinite cyclic) covering. Let
$\hat  
\alpha$ be a lift of
$\alpha$ to
$\hat  D_0$. Let $\hat  \beta $ be the lift of $\beta$ to $\hat  D_0$  which starts at $\hat  \beta (1)$.  Consider  a small neighborhood
$V\subset D $ of  the short arc  in $N$ connecting
$z_i$ to
$z'_i$ such that $V$ meets  $\alpha \beta$ only at $z'_i$. Let  $\gamma$ be a  generic loop  in $V\backslash z_i$
based at $z'_i$ which winds $w$ times around $z_i$ in the clockwise direction.  Let
$\hat 
\gamma$ be the  lift of
$\gamma$ to 
$\hat   D_0$ beginning at   $\hat  \beta(1)$ and ending at  $\hat  \alpha(0)$. 
We can assume that $\hat  \gamma$ is an embedded arc meeting $\hat  \alpha \hat  \beta$ only at
the  endpoints.
Let $\hat  z'_k$ be the first point of $\hat   \alpha$ which
lies also on  $\hat  \beta$ (possibly $\hat  z'_k= \hat  \alpha(1)$). Then $p(\hat  z'_k)=z'_k$ for some
$k=1,...,n$. Let
$\hat  \alpha'$ be the initial segment of $\hat  \alpha$ ending at $\hat  z'_k$. Let $\hat  \beta'$  be the final segment of $\hat  \beta$
starting at $\hat  z'_k$. Set $\hat  \delta= \hat  \alpha' \hat  \beta' \hat  \gamma$.  It is clear that  $\hat  \delta$ is a simple
closed curve in $\hat  D_0$.  By the Jordan curve theorem  it bounds a disc,  $  B\subset \hat  D_0$, which   lies either on the
left or on the right of $\hat 
\gamma$. Since
$\gamma$ passes clockwise around
$z_i$, 
  the  component of $D_0\backslash \Im(\gamma)$
adjacent to  $z_i$ lies on the right of $\gamma$. The set $p(B)\subset D_0$ being compact can not 
pass non-trivially over this component. Hence $  B$   lies on the left of $\hat  \gamma$. Therefore  $\hat  \delta$ passes
counterclockwise around $B$.

The number
$a_k-a_i$ is  equal to the
total winding number of the loop 
$p(\hat  \alpha') p(\hat  \beta') $ in  $D_n$ around  the
punctures $x_1,...,x_n$. Since  $\gamma$ is contractible in $D_n$,  
the loop 
$p(\hat  \alpha') p(\hat  \beta') $  is homotopic  in $D_n$ to
$p(\hat  \alpha') p(\hat  \beta')\gamma=p(\hat
\delta)$. Therefore $a_k-a_i$  is equal to the number of points in
$  B\cap p^{-1}(X)$. Since
$a_i$ is maximal, we must have 
$a_k=a_i$ and  
 $ B\cap p^{-1}(X)=\emptyset$, so that  $p(B)\subset  D_n\backslash z_i$. Then  we can isotop $T$ so as to have fewer
points of intersection with
$N$.  To see this we shall construct a bigon for the pair $(N,T)$. If $B$ meets $p^{-1} (N\cup T)$ only along  $\hat  \alpha'
\hat 
\beta'$ then the  projection
$p\vert_{\Int(B)}:
 \Int (B) \to  D_n\backslash z_i$ must be  injective. It follows that $w=1$ and  the union of $p(B)$ with the small disc bounded
by
$\gamma$ in $V$ is a bigon for   $(N,T)$.  Assume that $\Int(B)\cap p^{-1} (N\cup T)\neq \emptyset$. Note that
$p^{-1} (N)$ (resp.   $p^{-1} (T)$) is an  embedded one-manifold in $ \hat  D_0$ whose components are non-trivially permuted by
any covering transformation. If 
$p^{-1} (N)$ intersects $\Int(B)$  then  this intersection   consists of a finite set of disjoint arcs with endpoints on $\hat
\alpha'$. At least one of this arcs bounds together with a subarc of $\hat \alpha'$ a disc, $B_0\subset B$, whose interior does
not meet  $p^{-1} (N)$.  If 
$p^{-1} (N)$ does not meet  $\Int(B)$ then we set $B_0=B$. Applying the same construction to the
intersection of $B_0$ with  $p^{-1} (T)$ we obtain a bigon $B_{00} \subset B_0$ for the pair  
$(p^{-1} (N),p^{-1} (T))$. The restriction  of $p$ to $B_{00}$ is injective and yields a bigon for  $(N,T)$.
   Hence the intersection $N\cap T$ is not minimal. This contradicts   our choice of  $N,T$.
Therefore, the assumption 
$b_{i,i}<b_{i,j}$ must have been false. So, $b_{i,i}=b_{i,j}$.   Similarly, $b_{j,j}=b_{i,j}$. This
completes the proof of $(\ast)$ and of  Lemma 2.5.

\vskip 0.5truecm
\par \noindent {\bf 2.6. Lemma.} -- {\it  If  a self-homeomorphism $h$ of $D_n$  represents an element of 
$\Ker(B_n\to \Aut  (H_2(\tilde C)))$ then for any  noodle $N$ and  any fork $F$, we have $ \langle N,h(F)
\rangle =\langle N,F \rangle$.}

\vskip 0.5truecm
\par \noindent {Proof.} --  Let  $\{U_i\subset \Int(D)\}_{i=1}^m$ be disjoint closed disc neighborhoods of the points 
 $\{x_i\}_{i=1}^m$,   respectively. Let   $U$ be the set of points 
$\{x,y\}\in C$ such that at least one of $x, y$ lies in $\cup_{i=1}^m U_i$.  
Let 
$\tilde U \subset  \tilde C$ be the preimage of $U$ under  the covering map $ \tilde C \to C$. Observe  that the surface $\tilde
\S_F$ is an open square  such that for  a sufficiently big concentric closed subsquare  $S\subset \tilde \S_F$ we have 
  $\tilde \S_F\backslash S\subset \tilde U$.  Hence $\tilde \S_F$
represents a relative homology class 
$[\tilde
\S_F]\in H_2(\tilde C,
\tilde U)$. The boundary homomorphism $H_2(\tilde C, \tilde U) \to H_1(\tilde U)$ maps 
$[\tilde \S_F]$ into    $[\partial S]\in H_1(\tilde U)$. A direct computation in $\pi_1(U)$  (see [Bi2]) shows that
$(q-1)^2(qt+1)[\partial S]=0$. Therefore
$(q-1)^2(qt+1) [\tilde \S_F]= j(v_F)$ where $j$ is the inclusion homomorphism 
$H_2(\tilde C)  \to  H_2(\tilde C, \tilde U)$ and  
 $v_F\in H_2(\tilde C)$. 
Deforming if necessary $N$, we can assume that $N\cap (\cup_i U_i)=\emptyset$. Then    $\tilde
\S_N\cap \tilde U=\emptyset$ and therefore
$$(q-1)^2(qt+1) \langle N,F \rangle = \sum_{a,b\in {\bf Z}} (q^a t^b \tilde \S_N \cdot v_F) \,q^a t^b
 $$
where $q^a t^b \tilde \S_N \cdot v_F$ is the (well-defined) algebraic intersection number between a properly  embedded
surface and a 2-dimensional homology class. (This number does not depend on the choice of  $v_F$ as above).

 Any self-homeomorphism $h$ of $D_n$ is isotopic to a 
self-homeomorphism   of $D_n$  preserving the set $U$.  Therefore $v_{h(F)}=\tilde
h_*(v_F)$.
  If $\tilde h_*=\id$,   then
$$q^a t^b \tilde \S_N \cdot v_F=q^a t^b \tilde \S_N \cdot \tilde h_*(v_F)=
q^a t^b \tilde \S_N \cdot v_{h(F)}.$$ This implies that 
$(q-1)^2(qt+1) \langle N,F \rangle=(q-1)^2(qt+1) \langle N,h(F) \rangle$ and therefore
$ \langle N,F \rangle=  \langle N,h(F) \rangle$.

\vskip 0.5truecm
 \par \noindent {\bf 2.7. Deduction of Theorem  2.2 from the lemmas.} -- 
We shall prove that  a self-homeomorphism $h$ of $D_n$ 
representing an element of 
$\Ker(B_n\to \Aut  (H_2(\tilde C)))$  is isotopic to the identity map rel $\partial D_n$.
We begin with the following assertion.

\vskip 0.5truecm

($\ast\ast$)  {\it  An 
embedded arc 
$T$ in
$D_n$ with endpoints in (distinct) punctures can be isotopped off a noodle $N$ if and
only if
$h(T)$ can be isotopped off $N$.}

\vskip 0.5truecm

Indeed, we can extend 
$T$     to  a fork $F$ so that $T$ is the tines of $F$. 
By Lemma  2.6, $\langle N,h(F) \rangle
=0$ if and only if $\langle N,F \rangle =0$. Now Lemma  2.5   implies ($\ast\ast$).

We shall apply ($\ast\ast$) to the following arcs and noodles. 
For $i=1,..., n-1$, denote by $T_i$ the embedded arc  $ [x_i,x_{i+1}] \subset (-1,+1)\subset D$
and denote by $N_i$ the $i$-th \lq\lq elementary" noodle obtained by rushing from $d$ towards $x_i$, encircling $x_i$ in the
clockwise direction and then moving straight to $d'$. It is clear that $T_i$ can be isotopped off $N_j$ if and only if
$j\neq i, i+1$. This and    ($\ast\ast$)  imply that $h$ induces the identity permutation on the punctures of
$D_n$.

Since $T_1$ is disjoint from $N_3$, we can isotop $h$ rel $\partial D_n$ so that $h(T_1)$ is disjoint from $N_3$.
Similarly, $h(T_1)$ can be made disjoint from $N_4$. As it was explained in Sect.\  2.2, this can be done by a sequence of
isotopies eliminating bigons for the pair $(N_4,h(T_1))$. 
Since   $N_4$ and $h(T_1)$ do not meet  $N_3$,  neither  do these bigons.  Hence our  isotopies do not create intersections
of $h(T_1)$ with $N_3$. Repeating this argument, we can assume that $h(T_1)$ is disjoint from all $N_i$ with $i=3,4,..., n-1$.
By applying one final isotopy we can make
$h(T_1)=T_1$. 
Applying the same procedure to $T_2$ we can ensure that $h(T_2)=T_2$ while keeping $h(T_1)=T_1$.
Continuing in this way, we can  assume that $h(T_i)=T_i$ for all $i=1,...,n-1$.  Such a homeomorphism $h$ is isotopic to a
$k$-th power ($k\in {\bf Z}$) of the Dehn twist about a circle in  $ \Int (D_n)$  going very closely to
$\partial D_n$. This Dehn twist acts on $H _2(\tilde C)$ by multiplication by $q^{2n} t^2$. Since by
assumption $h$ acts trivially on $H _2(\tilde C)$, we must have $k=0$ so that $h$ is isotopic to the
identity rel $\partial D_n$.

 \vskip 0.5truecm
 \par \noindent {\bf 2.8.  Remarks.} -- The proof of Lemma 2.6 shows that each fork $F$
determines (a priori non-uniquely) a certain homology class  $v_F\in H_2(\tilde C)$. It follows from the computations in
[Bi2] that this class   is in fact well-determined by $F$. Thus, the forks yield a nice
geometric way of representing elements of   $H_2(\tilde C)$ (this was implicit  in [Kr1]).  For instance, for any $1\leq
i<j\leq n$ we can consider the fork  consisting of three linear segments connecting the point $-\sqrt {-1}/2 $ to $d,x_i,x_j$.
The corresponding classes
$\{v_{i,j}\in H_2(\tilde C) \}_{i,j}$ form a basis of the free $R$-module $ H_2(\tilde C)$. The action of the braid
generators $\sigma_1,...,
\sigma_{n-1}$ on this basis can be described by explicit formulas (see [Bi2], cf.  Sect.\  3.1).

\vskip 1.2truecm
\noindent {\bf 3. THE WORK OF   KRAMMER}
\vskip 0.5truecm
\par \noindent {\bf 3.1. A  representation of $B_n$.}  
Following [Kr2], we  denote by $ \Ref=\Ref_n$ the set of pairs of integers $(i,j)$ such
that $1\leq i<j\leq n$. Clearly, $\card (\Ref)= n (n-1)/2$. 

Let $R$ be a commutative ring with unit and $q,t\in
R$ be two invertible elements. 
Let $V=\oplus_{s\in \Ref} \,Rx_s$ be the free $R$-module  of rank $ n (n-1)/2$ with basis 
$\{x_{s}\}_{s\in \Ref}$.      Krammer [Kr2] defines
an
$R$-linear action of $B_n$ on $V$ by  
$$
\sigma_k(x_{i,j})=\cases {x_{i,j} & {{if}}   $\,k<i-1$ or $j<k$, \cr
x_{i-1,j}+  (1-q) x_{i,j} & {{if}}  $\,k=i-1$, \cr
tq(q-1) x_{i,i+1} + q x_{i+1,j} & {{if}}  $\,k=i <j-1$, \cr
tq^2 x_{i,j}  & {{if}}  $\,k=i =j-1$, \cr
x_{i,j} + tq^{k-i}  (q-1)^2 x_{k,k+1}& {{if}}  $\,i<k<j-1$, \cr
x_{i,j-1}+ t q^{j-i} (q-1) x_{j-1,j} & {{if}}  $\,k=j-1$, \cr
(1-q) x_{i,j} +q x_{i,j+1} & {{if}}  $\,k=j$, \cr}
$$ 
where $1\leq i<j\leq n$ and $k=1,...,  n-1$.  That the action of  $\sigma_k$ is invertible
and that 
relations  (0.1),  (0.2) are satisfied  should be verified by a direct computation.
For $R={{\bf Z}}[q^{\pm 1}, t^{\pm 1}]$, this  representation  is equivalent to the one considered in Sect.\  2.
In terms of the basis 
$\{v_{i,j}\in H_2(\tilde C) \}_{i,j}$ mentioned
in Sect.\  2.8, the equivalence   is given  by  
$$v_{i,j}=x_{i,j}+ (1-q) \sum_{i<k<j} x_{k,j},\,\,\, \,\,\, 
x_{i,j}=v_{i,j}+ (q-1) \sum_{i<k<j} q^{k-1-i} \,v_{k,j}.$$

\vskip 0.5truecm
\par \noindent {\bf 3.2. Theorem.} (D. Krammer   [Kr2])  -- {\it  
Let $R={\bf R}  [t^{\pm 1}]$, $q\in {\bf R}$, and $0<q<1$. Then the representation $B_n\to \Aut (V)$
defined in Sect.\ 3.1  is faithful for all $n\geq 1$.}

\vskip 0.5truecm

This Theorem implies Theorem 2.2:   if a representation over 
${{\bf Z}}[q^{\pm 1}, t^{\pm 1}]$ becomes faithful  after   assigning a real value to $q$, then it is faithful
itself. 

Below we outline  the main ideas of Krammer's proof. 

\vskip 0.5truecm
\par \noindent {\bf 3.3.   Positive braids and the set   $\Omega\subset B_n$.}  
We recall  
  a few    facts
about the braid group $B_n$, see [Ch], [Ga], [Mi]. 
For  $i=1,2,...,n-1$ denote by $s_i$ the  transposition  $(i,i+1) \in S_n$.
The set 
$\{s_1,...,s_{n-1}\}   $ generates the symmetric group $S_n$.  Let $\vert \, \, \vert :S_n\to {\bf  Z}$
be the length function with respect to this generating set: for $x\in S_n$,
$\vert  x \vert $ is the smallest natural number $k$ such that $x$ is a product of $k$ elements of the set 
$\{s_1,...,s_{n-1}\}$.  
The canonical projection $B_n\to S_n$ has a  unique  set-theoretic section
$r:S_n\to B_n$ such that $r(s_i)=\sigma_i$ for $i=1,...,n-1$
and $r(xy)= r(x)\, r(y)$ whenever $\vert  xy\vert =\vert  x\vert +\vert  y\vert$. The group $B_n$ admits a
presentation by generators
$\{ r(x) \, \vert\,  x\in S_n\}$ and relations $r(xy)= r(x) \,r(y)$  for all $x,y\in S_n$ such that
 $\vert  xy\vert =\vert  x\vert +\vert  y\vert$. Set $\Omega=  r(  S_n) \subset B_n$. Note that $1=r(1)\in
\Omega$ and $\sigma_i=r(s_i)\in
\Omega$ for all $i$.

The {\it positive braid monoid}  $B_n^+$ is the submonoid of $B_n$ generated by
$\{\sigma_1,...,\sigma_{n-1}\}$.  The elements of $B_n^+$ are called {\it
positive braids}.   Clearly, $\Omega\subset B_n^+$. 

For any $x\in B_n^+$ there is a unique longest $x'\in S_n$ such that  $x\in r(x') B_n^+$. We 
denote $r(x')\in \Omega$ by  $LF(x)$ where $LF$ stands for the leftmost factor. 
Observe that
$$LF(xy)=LF(x\, LF(y)) \eqno (3.1) $$
for any  $x,y\in  B_n^+$. This  implies that the map $B_n^+\times \Omega \to \Omega$ defined by
$(x,y)\mapsto LF(xy)$ is an action of the monoid $B_n^+$ on $\Omega$.

\vskip 0.5truecm
\par \noindent {\bf 3.4.   Half-permutations.}  
A set $A\subset \Ref$ is called a {\it half-permutation} if  whenever $(i,j), (j,k)\in A$, one has $(i,k)\in A$.
Each half-permutation $A$ determines
  an ordering $<_A$ on the set
$\{1,2,...,n\}$  by $i<_A j \Leftrightarrow (i,j)\in A$   (and vice versa).

To state deeper properties of half-permutations we  consider the set $2^{\Ref}$  of all subsets of
$\Ref$ and define a map
 $L: S_n\to 2^{\Ref}$ by
$$L(x)=\{ (i,j)\,\vert \, \, 1\leq
i<j\leq n, \,\,x^{-1}(i)> x^{-1}(j)\}\subset \Ref.$$
Note that  the set $  L(x)$ is a half-permutation and  $\card (L(x))=\vert x \vert$. It is obvious that   the
map 
$L: S_n\to 2^{\Ref}$ is injective.

  The key property of
half-permutations is the following assertion ([Kr2, Lemma 4.3]): for every half-permutation $A\subset \Ref$
there is a greatest  set $A'\subset A$ (with respect to inclusion) such that $A' =L(x)$ for a certain $x\in
S_n$. The corresponding  braid $r(x)\in \Omega$ is denoted by $GB(A)$ where $GB$ stands for the
greatest braid. This defines a map $GB$ from the set of half-permutations to $\Omega$.
In particular, for any $x\in S_n$ we have  $GB(L(x))=r(x)$.

\vskip 0.5truecm
\par \noindent {\bf 3.5.   Actions of $B_n$ of $ 2^{\Ref}$ and on half-permutations.}  
Let  $R={\bf R}  [t^{\pm 1}]$, $q\in {\bf R}$, and $0<q<1$.
The action of $B_n$  on $V$ defined in Sect.\  3.1  has the following property: for any  positive braid $x\in
B_n^+$ the entries of the  matrix  of the map $v\mapsto xv:V\to V$ belong to ${\bf R}_{\geq  0}+t{\bf
R}[t]$. (This is obvious for the generators $\sigma_1,...,\sigma_{n-1}$ of $B_n^+$). 
Therefore the action of $B_n^+$  preserves the set
$$\W=\bigoplus_{s\in \Ref}  ({\bf R}_{\geq  0}+t {\bf R}[t]) \,x_s \subset V.$$

 For a  set $A\subset \Ref$, define $\W_A$ as the subset of $\W$ consisting of  
vectors $\sum_{s\in \Ref}  k_s x_s$ with $k_s\in {\bf R}_{\geq  0}+t {\bf R}[t]$ such that 
$k_s\in  t {\bf R}[t]  \Leftrightarrow s \in A$. Clearly,  $\W$ is the disjoint union of
the sets $\W_A$ corresponding to various $A\subset \Ref$.
For any $x\in B_n^+$ and  $ A\subset \Ref$ there is a unique $B\subset \Ref$ such that $x\W_A\subset
\W_B$. We denote this set $B$ by $xA$. This defines an action of   $B_n^+$ on $2^{\Ref}$. By  [Kr2,
Lemma 4.2],  this action maps half-permutations  to half-permutations. This
defines an action of 
$B_n^+$ on the set of half-permutations.  Finally, Krammer observes that the map $GB$ from 
this set  to $\Omega$ is  $B_n^+$-equivariant. Thus,  for any positive braid $x\in
B_n^+$ and a half-permutation
$A\subset
\Ref$, we have
$$GB(xA)=LF(x\, GB(A)). \eqno (3.2)$$

The rest of the argument  is contained in the following two lemmas.

\vskip 0.5truecm
\par \noindent {\bf 3.6. Lemma.} -- {\it    Let $B_n$ act on a set $U$. Suppose we are given
non-empty disjoint   sets $\{ C_x \subset U\}_{x\in \Omega}$ such that $x C_y\subset 
C_{LF(xy)}$ for all $x,y\in \Omega$. Then the action of $B_n$ on $U$ is faithful. }

\vskip 0.5truecm
\par \noindent {Proof.} --  
Denote by $\ell$ the group homomorphism $B_n\to \Z$ mapping $\sigma_1,..., \sigma_{n-1}$ to $1$.
  We check first that the inclusion  $x C_y\subset 
C_{LF(xy)}$ holds for all $x\in B_n^+$ and  $y\in \Omega$.  Clearly, $\ell(B_n^+) \geq 0$ so that we can
use induction on $\ell(x)$. If $\ell(x)=0$ then  $x=1$  and the inclusion   follows from the  equality
$LF(y)=y$. Let
$\ell(x)\geq 1$. Then $x=\sigma_i v$ where 
$i=1,...,n-1$ and $v\in B_n^+$. Clearly, 
$\ell (v) =\ell(x)-1$. We have 
$$xC_y=\sigma_i vC_y\subset \sigma_i(C_{LF(vy)}) \subset C_{LF(\sigma_i LF(vy))}=C_{LF(\sigma_i
vy)}=C_{LF(xy)}.$$ Here the  first inclusion  follows from the induction hypothesis, the second 
inclusion  follows from the assumptions of the lemma, the middle equality
follows from (3.1).

It is known that for any $b\in B_n$ there are $x,y\in B_n^+$ such that 
$b=xy^{-1}$.  Therefore to prove the lemma it suffices to show that, if two elements $x,y\in B_n^+$
act in the same way on $U$, then $x=y$.  We will show this by induction on $\ell(x)+\ell (y)$. 
If $\ell(x)+\ell (y)=0$, then $x=y=1$.  Assume that $\ell(x)+\ell (y)>0$. 
By assumption, $C_1$ is non-empty; choose any $u\in C_1$. 
We have $xu\in xC_1 \subset C_{LF(x)}$ and similarly $yu\in  C_{LF(y)}$.
Hence $xu=yu\in  C_{LF(x)}\cap  C_{LF(y)}$. By the disjointness assumption, this is possible only if 
$LF(x)=LF(y)$.  Write $z=LF(x)\in \Omega $ and consider $x',y'\in B_n^+$ such that $x=zx'$ and
$y=zy'$.   We have $z\neq 1$ since otherwise $x=y=1$. Then $\ell(z)>0$ and   $\ell (x')+\ell (y')<\ell
(x)+\ell (y)$. The  induction  assumption yields that $x'=y'$. Therefore $x=y$. This proves the inductive
step and the lemma.

\vskip 0.5truecm
\par \noindent {\bf 3.7. Lemma.} -- {\it  For $ x\in \Omega$, set
$$C_x= \bigcup_{A \in GB^{-1} (x)}  \W_A \, \, \subset \,V$$ 
where $A$ runs over all half-permutations such that $GB(A)=x$. Then the sets $\{C_x\}_{x\in \Omega}$
satisfy all the conditions of Lemma 3.6. Therefore the action of $B_n$ on $V$ is faithful. } 

\vskip 0.5truecm
\par \noindent {Proof.} -- It is obvious that  the sets $\{C_x\}_{x\in \Omega}$ are disjoint. 
The set $C_x$ is non-empty because $\emptyset \neq \W_{L (r^{-1}(x))} \subset C_x$. 

To prove the inclusion $x C_y\subset 
C_{LF(xy)}$ for   $x,y\in \Omega$, it suffices to  prove that $x \W_A\subset 
C_{LF(xy)}$ whenever $A$ is a half-permutation  such that $GB(A)=y$.   
This   follows from the inclusions
$$x\W_A \subset \W_{xA}\subset C_{GB(xA)}= C_{LF(xy)}.$$
Here the first inclusion follows from the definition of the  $B_n^+$-action of $\Ref$. The second
inclusion follows from the definition of $C_{GB(xA)}$.  The  equality follows from 
(3.2). 

 \vskip 0.5truecm
 \par \noindent {\bf 3.8.  More about the representation.} -- 
Let $\ell_\Omega:B_n\to \Z$ 
be the length function with respect to the generating set $\Omega$: for $x\in B_n$, $  
\ell_\Omega (x)$ is the minimal natural number $k$ such that $x=x_1...x_k$ where
 $x_i \in 
 \Omega$ or $x_i^{-1}\in \Omega$ for each $i=1,...,k$.  
Among other related results, Krammer  gives an explicit  computation of $\ell_\Omega  $
in terms of his representation. Namely, 
 take the Laurent polynomial ring  $R={{\bf Z}}[q^{\pm 1}, t^{\pm 1}]$ as the ground ring
and denote by $\rho$ Krammer's representation $B_n\to \Aut(V)$ defined in Sect.\  3.1.
For $x\in B_n$, consider the Laurent expansion $\rho(x)= A_kt^k+A_{k+1} t^{k+1}+...+A_l t^l$
where $\{A_i\}_{i=k}^l$ are $(m\times m)$-matrices over ${\bf Z} [q^{\pm 1}]$ and $A_k\neq 0, A_l\neq 0$.
Then
$$\ell_\Omega (x)= \max\, (l-k,l,-k).$$
This formula   yields another proof of the faithfulness of $\rho$: If $x\in \Ker \rho$, then
$k=l=0$ so that $\ell_\Omega(x)=0$ and  $x=1$.

The length function   $  
\ell_\Omega $ was first considered by Charney [Ch] who proved that the formal power series  
$\sum_{x\in B_n} z^{\ell_\Omega (x)} \in {\bf Z}[[z]]$
is a rational function. It is unknown whether   the similar formal power series  determined by
the length function with respect to the generators $\sigma_1,..., \sigma_{n-1}$ is rational.

\vskip 1.2truecm
\noindent {\bf  4.    BMW-ALGEBRAS AND   REPRESENTATIONS OF $B_n$}
\vskip 0.5truecm
\par \noindent {\bf 4.1. Hecke algebras and representations of $B_n$.}  A vast family of
representations of $B_n$ including the Burau representation arise from a study of Hecke algebras. The
Hecke algebra
$H_n (\alpha)$ corresponding to  
$\alpha\in {\bf C}$  can be defined  as  the quotient of the complex group ring ${\bf C} [B_n]$ by the
relations
$\sigma_i^2= (1-\alpha) \sigma_i + \alpha$  for $i=1,..., n-1$. This  family of
 finite dimensional  ${\bf C}$-algebras    is a
one-parameter deformation of    ${\bf
C} [S_n]=H_n (1)$.  For $\alpha$ sufficiently close to $1$,  the algebra
$H_n (\alpha)$ is isomorphic to
${\bf
C} [S_n]$. For such $\alpha$, the algebra  $H_n (\alpha)$ is semisimple and its irreducible
representations are indexed by  the Young diagrams  with $n$ boxes. Their decomposition rules and
dimensions are the same  as for  the  irreducible representations of $S_n$.  

Each   representation of
  $H_n (\alpha)$ yields a
representation of
$B_n$ via the natural projection $B_n\subset {\bf C} [B_n] \to  H_n (\alpha)$. This gives  a family of 
irreducible finite dimensional  representations of $B_n$ indexed by the Young diagrams with $n$ boxes.
These representations were extensively studied by V. Jones [Jo].  In particular, he observed that the 
$(n-1)$-dimensional Burau representation of
$B_n$    appears as the irreducible representation associated
with the two-column Young diagram whose  columns have  $n-1$ boxes and   one box, respectively.

\vskip 0.5truecm
\par \noindent {\bf 4.2.  Birman-Murakami-Wenzl algebras and their representations.}  
Jun Murakami [Mu] and independently  J.  Birman
and H. Wenzl [BW]  introduced a   two-parameter family of  finite dimensional    ${\bf C}$-algebras $C_n
(\alpha, l)$  where  $\alpha$ and $l$ are 
 non-zero complex
numbers such that $\alpha^4\neq 1$ and $l^4\neq 1$. For
$i=1,...,n-1$, set 
$$e_i= (\alpha+\alpha^{-1})^{-1} (\sigma_i+\sigma_i^{-1})-1\in {\bf
C} [B_n].$$ 
The algebra $C_n (\alpha, l)$ is the quotient of $ {\bf
C} [B_n]$ by the relations 
$$e_i \sigma_i =l^{-1} e_i,\,\,\, \,\,\,e_i \sigma_{i-1}^{\pm 1}  e_i =l^{\pm 1} e_i$$
for all  $i$.  (The original definition in [BW]  involves more relations; for  the shorter list given
above, see [We]). The algebra $C_n (\alpha, l)$ admits a geometric interpretation in terms of so-called 
Kauffman skein classes of tangles in Euclidean 3-space. This family of algebras is a deformation of an
algebra introduced by R. Brauer [Br] in 1937. 

The algebraic structure and  representations of $C_n (\alpha, l)$ were  studied by  Wenzl
[We], who established (among other results) the following three facts. 

\vskip 0.5truecm

(i) {\it For generic $\alpha, l$, the algebra
$C_n (\alpha, l)$ is semisimple.}

\vskip 0.5truecm

  Here
\lq\lq generic" means that  $\alpha$ is not a root of unity and $\sqrt {-1}\,l $ is not an integer power of
$- \sqrt {-1} \,\alpha$. (The latter two numbers  correspond to $r$ and $q$ in Wenzl's notation). 
In the sequel
we   assume that
$\alpha, l$ are generic  in this sense. We   denote the number of
boxes in  a Young diagram
$\lambda$ by
$\vert
\lambda
\vert$.  

\vskip 0.5truecm

(ii) {\it  The  irreducible  finite dimensional   $C_n (\alpha, l)$-modules 
are indexed by  the Young diagrams  $\lambda$ such that   $\vert \lambda \vert \leq n$ and $\vert \lambda
\vert\equiv n \,(\mod 2)$.  }

\vskip 0.5truecm

The irreducible
$C_n (\alpha, l)$-module corresponding to   $\lambda$ will be denoted by 
$V_{n,\lambda}$. Composing the natural projection $B_n\subset {\bf C} [B_n] \to  C_n (\alpha,l)$ with  the
action of
$C_n (\alpha,l)$ on $V_{n,\lambda}$
we obtain an irreducible representation of
$B_n$.

Observe that the  inclusion $B_{n-1} \hookrightarrow
B_{n}$ sending each
$\sigma_i
\in B_{n-1}$ with $i=1,...$, $n-2$ to
$\sigma_i \in B_{n }$ induces an inclusion $C_{n-1} (\alpha, l) \hookrightarrow  C_{n} (\alpha, l)$ for all
$n\geq 2$.

\vskip 0.5truecm

(iii) {\it The $C_n (\alpha, l)$-module  $V_{n,\lambda}$
  decomposes as a   $C_{n-1} (\alpha, l)$-module  into a direct sum 
  $\oplus_{\mu} V_{n-1,\mu}$ where $\mu$ ranges over all Young
diagrams obtained by removing or (if $\vert \lambda \vert<n$)  adding  one box to $\lambda$. 
Each such   $\mu$  appears in this decomposition with multiplicity 1. }

\vskip 0.5truecm
\par \noindent {\bf 4.3.  The Bratelli diagram  for the BMW-algebras.}  
The  assertions (ii) and (iii) in Sect.\  4.2 allow  us to draw the Bratelli diagram  for the sequence  
$C_1 (\alpha, l)\subset  C_2 (\alpha, l) \subset ...$  
On the level $n=1,2,...$ of the Bratelli diagram   one puts all  Young diagrams 
 $\lambda$ such that   $\vert \lambda \vert \leq n$ and $\vert \lambda
\vert\equiv n \,(\mod 2)$.  Then one connects  by an edge each $\lambda$ on the $n$-th level to all 
Young diagrams
on the $(n-1)$-th  level    
obtained by removing or (if $\vert \lambda \vert<n$)  adding  one box to $\lambda$. 
For instance, the $n=1$ level consists  of  a single Young diagram with one box corresponding to the
tautological
 one-dimensional representation of  $C_{1} (\alpha, l)={\bf C}$. The $n=2$ level contains the empty Young
diagram and two Young diagrams with two boxes. All three are connected  by an edge  to the   diagram on
the 
level 1. 
 Note that every Young diagram $\lambda$ appears on the levels   $\vert \lambda \vert, \vert \lambda
\vert +2, \vert \lambda \vert+4,...$ 

 The Bratelli diagram  yields a useful method of computing the dimension  of  
  $V_{n,\lambda}$ where $\lambda$ is a Young diagram on the $n$-th level.  It is clear
from (iii) that
 $\dim (V_{n,\lambda}) $  is the number of paths on the Bratelli diagram  leading  from $\lambda$ to the
only diagram on  the  level 1. Here by a path we mean a path with   vertices lying on consecutively
decreasing levels. We give  three examples of computations based on (iii). 

(a) Let  $\lambda_n$  
be 
the Young
diagram    represented by a column of   $n$ boxes. There is only one path from
$\lambda_n$,  
 positioned on the level $n$,   to the top of  the Bratelli diagram. Hence, $\dim (V_{n,\lambda_n}) =1$ for all
$n\geq 1$. 

For $n\geq 2$, the algebra   $C_{n} (\alpha, l)$   has two
one-dimensional representations. In both of them all $e_i$ act as $0$ and  all $\sigma_i$ act as
multiplication by one and  the same number equal  either to $\alpha$ or  to $\alpha^{-1}$. We choose the
correspondence between the irreducible $C_{n} (\alpha, l)$-modules and the Young diagrams so that
all $\sigma_i$ act on $V_{n,\lambda_n}$ as multiplication by $\alpha$.  If $\lambda^T_n$ is the
  Young diagram     represented by a   row of   $n$ boxes,  then similarly to (a) we have that $\dim (V_{n,
\lambda^T_n}) =1$   and all
$\sigma_i$ act on $ V_{n, \lambda^T_n}$ as multiplication by $\alpha^{-1}$.

(b) For $n\geq 2$, let  $\lambda'_n$
be the  two-column Young diagram  whose   columns have $n-1$ boxes and one box, respectively.
For $n\geq 3$, the diagram  $\lambda'_n$,  positioned on the level $n$,  is
connected to only two Young diagrams on the previous level, namely, to $\lambda'_{n-1}$ and
$\lambda_{n-1}$. Hence $$ \dim (V_{n, \lambda'_n})=\dim (V_{n-1, \lambda'_{n-1}})+\dim (V_{n-1,
\lambda_{n-1}}) =\dim (V_{n-1, \lambda'_{n-1}})+1.$$
We have $\lambda'_{2}=\lambda^T_{2}$ so that  $\dim (V_{2, \lambda'_{2}})=1$.
Hence  $ \dim (V_{n,
\lambda'_n})=n-1$ for all
$n\geq 2$.

(c)  For $n\geq 2$, consider the module $V_{n, \lambda_{n-2}}$ corresponding to the Young diagram
 $\lambda_{n-2}$ positioned on the level $n$. 
 If $n\geq 3$, then this diagram    is
connected to three Young diagrams on the previous level, namely, to $\lambda_{n-1}, \lambda_{n-3}$, 
$\lambda'_{n-1}$. Hence $$ \dim (V_{n, \lambda_{n-2}})=\dim  (V_{n-1, \lambda_{n-1}})+
\dim (V_{n-1, \lambda_{n-3}})+\dim (V_{n-1,
\lambda'_{n-1}})$$
$$ =\dim (V_{n-1, \lambda_{n-3}})+n-1.$$
We gave $\lambda_0=\emptyset $ and by (iii) above,   
  $\dim (V_{2, \lambda_{0}})=\dim (V_{1, \lambda_{1}})=1$. Thus for all $n\geq 2$,  $$ \dim (V_{n,
\lambda_{n-2}})=n(n-1)/2,$$  i.e.,   $ V_{n, \lambda_{n-2}}$ has the same  dimension 
  as the Krammer representation of $B_n$. We now rescale the 
representation   $B_n\to \Aut  (V_{n,
\lambda_{n-2}})$
by  dividing the action of each $\sigma_i$ by $\alpha$.

\vskip 0.5truecm
\par \noindent {\bf 4.4. Theorem.} (M. Zinno [Zi]) -- {\it  The 
Krammer representation corresponding to $q=-\alpha^{-2}$ and $ t=\alpha^3
l^{-1}$ is isomorphic  to the rescaled
representation $B_n\to \Aut  (V_{n,
\lambda_{n-2}} )$. } 
\vskip 0.5truecm

 The proof given in  [Zi] goes by a direct comparison of both actions of $B_n$ on  certain bases.
Theorem 4.4 implies that the Krammer-Bigelow representation considered in Sect.\  2 and 3 
is irreducible.

\vskip 1.5truecm
\centerline{\bf BIBLIOGRAPHY}
\vskip 0.5truecm

\item{[Bi1]} S. BIGELOW -- {\it 
The Burau representation   is not faithful for  $n\geq 5$},
Preprint math. GT/9904100.

\item{[Bi2]} S. BIGELOW -- {\it 
Braid groups are linear},
Preprint  math. GR/0005038.

\item{[Bir]} J. S.  BIRMAN -- {\it Braids, links, and mapping class   groups},
Ann. of Math.  Stud.,  vol. 82,
Princeton Univ. Press, Princeton, N.J., 1974.

\item{[BW]} J. S.  BIRMAN, H. WENZL -- {\it Braids, link polynomials and a new algebra},
Trans. Amer. Math. Soc. {\bf 313} (1989),   249--273.

\item{[Br]}  R. BRAUER -- {\it On algebras which are connected with the semisimple continuous groups}, 
Ann. of Math.    
{\bf 38} (1937),  
  857--872.

\item{[Bu]}  W. BURAU -- {\it  \"Uber Zopfgruppen und gleichsinnig verdrillte Verkettungen},
Abh. Math. Semin. Hamburg. Univ. {\bf 11} (1935),  
179-186.

\item{[Ch]} R. CHARNEY -- {\it  Geodesic automation and growth functions for Artin groups of finite type},
Math. Ann. {\bf 301} (1995), 307--324.

\item{[DFG]}  J. L. DYER, E. FORMANEK, E. K. GROSSMAN -- {\it  On the linearity of automorphism groups of free groups},
Arch. Math. (Basel) {\bf  38} (1982),  
404--409.

\item{[FLP]} A. FATHI, F. LAUDENBACH, V. POENARU -- {\it 
Travaux de Thurston sur les surfaces}, Ast\'erisque
{\bf  66--67}
(1991), Soc. Math. France, Paris.

\item{[FP]}    E. FORMANEK, C. PROCESI -- {\it  The automorphism group of a free group is not linear},
J. Algebra  {\bf  149} (1992),  
494--499.

\item{[Ga]} F. A. GARSIDE -- {\it The braid group and other groups}, Quart. J. Math. Oxford {\bf 20}
(1969), 235--254.

\item{[Iv1]} N. V. IVANOV -- {\it Automorphisms of Teichm\"uller modular groups}, 
Topology
and geometry -- Rohlin Seminar,  Lecture Notes in Math.,  vol. 1346,
Springer, Berlin-New York (1988), 199--270.

\item{[Iv2]} N. V. IVANOV -- {\it  Subgroups of Teichm\"uller modular groups},
Transl.\ of Math.  Monographs, vol. 115.
Amer. Math. Soc., Providence, RI, 1992. 

\item{[Jo]} V. F. R. JONES -- {\it   Hecke algebra representations of braid groups and link polynomials},
Ann. of Math.   {\bf 126}  (1987),   335--388. 

\item{[Ka]}  C. KASSEL -- {\it L'ordre de Dehornoy sur les tresses}, 
S\'eminaire Bourbaki, expos\'e 865  (novembre 1999), to appear in Ast\'erisque, Soc. Math. France, Paris.

\item{[Kr1]} D.  KRAMMER -- {\it 
The braid group $B_4$ is linear},
Preprint, Basel  (1999).

\item{[Kr2]} D.  KRAMMER -- {\it 
Braid Groups are linear},
Preprint, Basel  (2000).

\item{[La]} R. J. LAWRENCE -- {\it Homological representations of the Hecke algebra}, Comm. Math. Phys.  {\bf 135} 
(1990),   141--191. 

\item{[LP]} D. D. LONG, M. PATON -- {\it 
The Burau representation is not faithful for $n\geq 6$}, Topology 
{\bf  32}
(1993), 439--447.

\item{[Mi]} J. MICHEL -- {\it A note on words  in braid monoids}, J. Algebra {\bf 215}
(1999), 366--377.

\item{[Mo]} J. A. MOODY -- {\it 
The Burau representation of the braid group $B\sb n$ is unfaithful for large $n$},
Bull. Amer.
Math. Soc. (N.S.)  
{\bf  25}
  (1991),   379--384.

\item{[Mi]} J. MURAKAMI -- {\it The Kauffman polynomial of links and representation theory}, 
Osaka J. Math.
{\bf 24} (1987),  745--758.

\item{[We]}  H. WENZL -- {\it Quantum groups and subfactors of type $B$, $C$, and $D$}, 
 Comm. Math. Phys. 
{\bf 133} (1990),  
383--432.

\item{[Zi]} M. G. ZINNO -- {\it 
On Krammer's Representation of the Braid Group},
Preprint math. RT/0002136.

\vskip 1truecm
\hskip 6.5truecm Vladimir TURAEV\par
\medskip
\hskip 6.5truecm Institut de Recherche Math\'ematique Avanc\'ee\par
\hskip 6.5truecm  Universit\'e Louis Pasteur et C.N.R.S.\par
\hskip 6.5truecm 7 rue Descartes\par
\hskip 6.5truecm F-67084 STRASBOURG Cedex\par
\smallskip
\hskip 6.5truecm E--mail~: turaev@math.u-strasbg.fr\par

\end

 One can check that a set $A\subset \Ref$ belongs to the image of $L$ if
and only if both $A$ and $\Ref \backslash A$ are half-permutations.

\item{[1]} N. AUTEUR - {\it titre de l'article}, revue {\bf num\'ero}
(ann\'ee), pages.

\item{[1]} L.V. AHLFORS -- {\it title.} McGraw-Hill, New York 1973

\bibitem{BL1} K. BURDZY, G.F. LAWLER -- {\it Title} Probab. Th. Rel. Fields
{\bf 84}
(1990), 393-410.

\item{[BLM]} J. S.  BIRMAN, D. D. LONG, J. A. MOODY  -- {\it Finite-dimensional
representations of Artin's braid group},   The mathematical legacy of Wilhelm Magnus: groups,
geometry and special functions (Brooklyn, NY, 1992), Contemp. Math., vol. 169, Amer.
Math. Soc., Providence, RI (1994), 123--132.

 $ where $[a,b]=a^{-1} b^{-1} ab]$. 

The linearity of $B_n$ implies that $B_n$ satisfies  all the properties of linear groups
including the Tits alternative,

Note that $\tau_\alpha^2=\tau_\beta$ where $\beta$ is the simple closed curve 
in $\S$ bounding a 
regular
neighborhood of $\alpha$.

$\alpha$ and $\beta$ are disjoint up to
isotopy, i.e.,

If $\alpha\subset \S$ is a simple closed curve bounding a compact subsurface of $\S$
then $[\alpha]=0$ and therefore $(\tau_\alpha)_*=\id_H$.
However,  $\tau_\alpha\neq 1$ unless $\alpha$  bounds a disc in $\S$. 

Let $\overline  C$ (resp. $C$) be the space of all  ordered  (resp. unordered) pairs of distinct points in $D_n$.
In other words, $\overline  C=(D_n\times D_n)\backslash \hbox{diagonal}$ and $C$ is obtained from    $\overline  C$ by the 
identification $(x,y)=(y,x)$ for  any $x, y \in D_n$ with $x\neq  y$. It is clear that $C$ and $\overline  C$ are connected
non-compact 4-manifolds with boundary. 
They  have a  natural orientation induced by the counterclockwise  orientation of
$D_n$.  Set $d=-i\in \partial D$, $d'= -i e^{ \delta { {\pi i} \over {2}}}\in \partial D$  with small
positive $\delta$  and take
$c_0=\{d,d'\}$ as the basepoint for $C$. Here $i=\sqrt{-1}\in {\bf C}$.

\vskip 0.5truecm
\par \noindent {Proof.} --  Let  $\{U_i\subset \Int(D)\}_{i=1}^m$ be disjoint closed disc neighborhoods of the points 
 $\{x_i\}_{i=1}^m$,   respectively. Let   $U$ be the set of points 
$\{x,y\}\in C$ such that at least one of $x, y$ lies in $\cup_{i=1}^m U_i$.  
Let 
$\tilde U \subset  \tilde C$ be the preimage of $U$ under  the covering map $ \tilde C \to C$. Observe  that the surface $\tilde
\S_F$ is an open square  such that for  a sufficiently big concentric closed subsquare  $S\subset \tilde \S_F$ we have 
  $\tilde \S_F\backslash S\subset \tilde U$.  Hence $\tilde \S_F$
represents a relative homology class 
$[\tilde
\S_F]\in H_2(\tilde C,
\tilde U)$. The boundary homomorphism $H_2(\tilde C, \tilde U) \to H_1(\tilde U)$ maps 
$[\tilde \S_F]$ into    $[\partial S]\in H_1(\tilde U)$. A direct computation in $\pi_1(U)$  (see [Bi2]) shows that
$(q-1)^2(qt+1)[\partial S]=0$. Therefore
$(q-1)^2(qt+1) [\tilde \S_F]= j([\tilde \S'_F])$ where $j$ is the inclusion homomorphism 
$H_2(\tilde C)  \to  H_2(\tilde C, \tilde U)$ and  
 $[\tilde \S'_F]$ is an element of $H_2(\tilde C)$
represented by  a closed oriented immersed surface  $\tilde \S'_F$ in $\tilde C$. 
deforming if necessary $N$ we can assume that $N\cap \cup_i U_i=\emptyset$. Then    $\tilde \S_N\cap \tilde
U=\emptyset$ and therefore
$$(q-1)^2(qt+1) \langle N,F \rangle = \sum_{a,b\in {\bf Z}} (q^a t^b \tilde \S_N \cdot \tilde \S'_F) q^a t^b
.$$
Here $q^a t^b \tilde \S_N \cdot \tilde \S'_F$ is the (well-defined) algebraic intersection number between a properly  embedded
surface and an immersed closed surface. (This number does not depend on the choice of  $\tilde \S'_F$ as above).
For any self-homeomorphism $h$ of $D_n$ we have that
$\tilde \S'_{h(F)}=\tilde h(\tilde \S'_F)$. If $:tilde h$ acts trivially on $H_2(\tilde C)$ then
$q^a t^b \tilde \S_N \cdot \tilde \S'_F=q^a t^b \tilde \S_N \cdot \tilde h(\tilde \S'_F)=
q^a t^b \tilde \S_N \cdot \tilde \S'_{h(F)}$. This implies that 
$(q-1)^2(qt+1) \langle N,F \rangle=(q-1)^2(qt+1) \langle N,h(F) \rangle$ and therefore
$ \langle N,F \rangle=  \langle N,h(F) \rangle$.

\vskip 0.5truecm
\par \noindent {\bf 3.4.   Construction of the sets  $\{ C_x \subset V\,\vert\, x\in \Omega\}$.}  We
consider now the representation $B_n\to \Aut (V)$ defined above and ouline  the construction of the sets 
$\{ C_x
\subset V\,\vert\, x\in
\Omega\}$ satisfying Lemma 3.3. Denote by $ \Ref$ the set of pairs of integers $(i,j)$ such that $1\leq
i<j\leq n$.  Define a map
 $L: S_n\to 2^{\Ref}$ by
$$L(x)=\{ (i,j)\,\vert \,  1\leq
i<j\leq n, x^{-1}(i)> x^{-1}(y)\}.$$
(Here $2^{\Ref}$ is the set of subsets of $\Ref$). 
Note that  the set $  L(x)$ has $\e(x)$ elements. It is easy to see that the map  $L: S_n\to 2^{\Ref}$ is
injective.

$$\rho(x)=\sum_{i=k}^l A_i  t^i$$

as
a subgroup of the automorphism group of a free group  of rank $n$,